\input amstex
\loadeufm

\documentstyle{amsppt}

\magnification=\magstep1

\baselineskip=20pt
\parskip=5.5pt
\hsize=6.5truein
\vsize=9truein
\NoBlackBoxes

\define\br{{\Bbb R}}

\define\e{{\varepsilon}}
\define\OO{{\Omega}}
\define\CL{{\Cal{L}}}

\topmatter
\title
 Necessary and Sufficient Conditions
for the Solvability of the $L^p$ Dirichlet Problem
On Lipschitz Domains. 
\endtitle

\author Zhongwei Shen
\endauthor

\leftheadtext{Zhongwei Shen}
\rightheadtext{The $L^p$ Dirichlet Problem}
\address Department of Mathematics, University of Kentucky,
Lexington, KY 40506.
\endaddress

\email shenz\@ms.uky.edu
\endemail

\abstract
We study the homogeneous elliptic systems 
of order $2\ell$ with real constant coefficients
on Lipschitz domains in $\br^n$, $n\ge 4$. 
For any fixed $p>2$, we show that a reverse H\"older condition with 
exponent $p$ is
necessary and sufficient 
for the solvability of the Dirichlet problem with boundary
data in $L^p$.
We also obtain a simple sufficient condition.
As a consequence, we establish the solvability of the 
$L^p$ Dirichlet problem
for $n\ge 4$ and $2-\e< p<\frac{2(n-1)}{n-3} +\e$.
The range of $p$ is known to be sharp if $\ell\ge 2$ and
$4\le n\le 2\ell +1$.
For the polyharmonic equation, 
the sharp range of $p$ is also found in the case
$n=6$, $7$ if $\ell=2$, and $n=2\ell+2$ 
if $\ell\ge 3$.
\endabstract



\endtopmatter

\document

\centerline{\bf 1. Introduction}

In this paper we study the higher order homogeneous
elliptic systems
with real constant coefficients
on bounded domains in $\br^n$ with Lipschitz boundaries. 
For any fixed $p>2$, we obtain necessary and sufficient conditions
for the solvability of the Dirichlet problem with boundary
data in $L^p$. As a consequence, we are able to establish
the solvability of the $L^p$ Dirichlet problem
for $2-\e<p<\frac{2(n-1)}{n-3}+\e$.
The range of $p$ is known to be sharp if
$\ell\ge 2$ and $4\le n\le 2\ell +1$, where $2\ell$ is the order of the
system. We also obtain the $L^p$ solvability for the sharp
range of $p$ in the case of the polyharmonic equation $\Delta^\ell u=0$
for $n=6$, $7$, if $\ell=2$ and  for $n=2\ell +2$,
if $\ell\ge 3$.

More precisely,
let $\OO$ be a bounded Lipschitz domain in $\br^n$.
Consider the homogeneous elliptic system of order $2\ell$,
$\CL(D) u =0$ in $\OO$, where $u=(u^1,\cdots, u^m)$,
$$
\aligned
\big(\CL (D) u\big)^j &=\sum_{k=1}^m \CL^{jk}(D) u^k,\ \ 
j=1,\dots, m,\\
\CL^{jk}(D)
&=\sum_{|\alpha|=|\beta|=\ell}
a_{\alpha\beta}^{jk} D^\alpha D^\beta,
\endaligned
\tag 1.1
$$
and $D=(D_1,D_2,\dots, D_n)$, $D_i=\partial/\partial x_i$
for $i=1,2,\dots,n$. Also $\alpha=(\alpha_1,\alpha_2,\dots,\alpha_n)$
is a multi-index with length $|\alpha|=\alpha_1+\cdots+\alpha_n$,
and $D^\alpha =D_1^{\alpha_1}D_2^{\alpha_2}\cdots D_n^{\alpha_n}$.
Let
$$
\CL^{jk} (\xi)=\sum_{|\alpha|=|\beta|=\ell}
a_{\alpha\beta}^{jk} \xi^{\alpha}\xi^\beta,
\ \ \ \text{ for }\xi\in \br^n.
$$
Throughout this paper, we will assume that $a_{\alpha\beta}^{jk}$
are real constants satisfying the symmetry condition
$$
\CL^{jk}(\xi)=\CL^{kj}(\xi)
\tag 1.2
$$
 and the Legendre-Hadamard 
ellipticity condition,
$$
\mu |\xi|^{2\ell} |\eta|^2\le \sum_{j,k=1}^{m}
\CL^{jk}(\xi)
\eta_j \eta_k
\le \frac{1}{\mu}\, |\xi|^{2\ell} |\eta|^2,
\tag 1.3
$$
for some $\mu>0$, and all $\xi\in \br^n$, $\eta\in\br^m$.
We are interested in the Dirichlet problem,
$$
\left\{
\aligned
\CL (D) u & =0 \ \ \ \ \ \text{ in } \OO,\\
D^\alpha u &=f_\alpha \ \ \ \text{ on }
\partial\OO \ \ \ \
\text{ for } |\alpha|\le \ell -2,\\
\frac{\partial^{\ell -1} u}{\partial N^{\ell -1}}
&=g \ \ \ \ \ \text{ on } \partial\OO,
\endaligned
\right.
\tag 1.4
$$
where $\frac{\partial^{\ell -1} u}{\partial N^{\ell -1}}
=\sum_{|\alpha|=\ell -1}
\frac{(\ell -1)!}{\alpha !} N^{\alpha} D^\alpha u$,
and $N$ denotes the outward unit normal to $\partial\OO$.

To describe the $L^p$ Dirichlet problem,
we let $\dot{f}=\{ f_\alpha:\ |\alpha|\le \ell -2\}$ be an array
of functions on $\partial\OO$. Following Verchota and Pipher
\cite{V3,PV3,V4}, we consider the Dirichlet problem (1.4) with 
 boundary data $(\dot{f}, g)$ taken from the
space $\big(WA^p_{\ell-1}(\partial\OO), L^p(\partial\OO)\big)$, 
where $WA^p_{\ell-1}(\partial\OO)$
denotes the completion of the set of arrays
 $\dot{\psi}=\{ D^\alpha \psi:  |\alpha|\le \ell-2\}$,
$\psi\in C^\infty_0 (\br^n)$ under the norm 
$$
\sum_{|\alpha|\le \ell-2} \|D^\alpha \psi\|_{L^p(\partial\OO)}
+\sum_{|\alpha|=\ell-2} \| \nabla_t D^\alpha \psi\|_{L^p(\partial\OO)}
\ \ \ \ \ \ \text{ on }\ \partial\OO.
\tag 1.5
$$
 Here $\nabla_t h$ denotes the tangential
derivatives of $h$. The boundary values in (1.4)
are taken in the sense of non-tangential convergence a.e.
with respect to the surface measure $d\sigma$ on $\partial\OO$.
As such, we will require that the non-tangential maximal
function $(\nabla^{\ell-1} u)^*$ on $\partial\OO$
is in $L^p(\partial\OO)$,
where $\nabla^{\ell-1} u$ denotes the tensor of
all partial derivatives of order $\ell-1$ in $\br^n$ 
of $u$. Thus the $L^p$ Dirichlet problem (1.4) on $\OO$ is said to be uniquely
solvable if given any $(\dot{f},g)\in \big(WA^p_{\ell-1}(\partial\OO),
L^p(\partial\OO)\big)$, there exists a unique $u$ satisfying (1.4)
and $(\nabla^{\ell-1} u)^*\in L^p(\partial\OO)$, and the
unique solution $u$ satisfies the scale-invariant estimate
$$
\| (\nabla^{\ell -1}u)^*\|_{L^p(\partial\OO)}
\le C\, \bigg\{
\| g\|_{L^p(\partial\OO)}
+\sum_{|\alpha|=\ell -2}
\|\nabla_t f_\alpha \|_{L^p(\partial\OO)}\bigg\}
\tag 1.6
$$
with constant $C$ independent of the boundary data $(\dot{f},g)$.

For Laplace's equation $\Delta u=0$ in $\OO$, the $L^p$ Dirichlet
problem was solved by Dahlberg \cite{D1,D2}
for the optimal range $2-\e<p\le \infty$, where $\e>0$
depends on $n$ and the Lipschitz character of $\OO$
(also see \cite{JK,DK1,K1} for the Neumann problem
and \cite{K2} for other related problems).
For the general elliptic equations
and systems $\CL(D) u=0$ considered in this paper,
 the solvability of the $L^p$ Dirichlet problem 
has been established
for $2-\e<p<2+\e$. See \cite{FKV,DKV2,K1,F,G} for second order 
elliptic systems, \cite{DKV1,V2,V3} for the biharmonic and polyharmonic
equations, and \cite{PV3,V4} for general higher order
elliptic equations and systems.
It is known that the restriction $p>2-\e$ is necessary for 
general Lipschitz domains (see e.g. \cite{K1}).
However, due to the lack of the maximum principle for
elliptic systems and higher order elliptic equations,
it has been a challenging problem
 to determine the optimal ranges of $p$
for which one may solve the $L^p$ Dirichlet problem.
Nevertheless
in the case of $n=2$ or $3$, the $L^p$ Dirichlet problem
for elliptic systems and higher order equations was solved
for the optimal range $2-\e<p\le \infty$ \cite{PV1,DK2,PV2,PV4,V4,S1,S2,MM}.
This was done by establishing certain decay
estimates on the Green's
functions. In the lower dimensional case, these estimates
lead to the Miranda-Agmon maximal principle
$\| \nabla^{\ell -1} u\|_{L^\infty(\OO)}
\le C\, \|\nabla^{\ell -1} u\|_{L^\infty(\partial\OO)}$,
from which the $L^p$ solvability for $2<p<\infty$ follows
by interpolation. However, in the case $n\ge 4$, these decay
estimates only yield the solvability of the Dirichlet
problem in certain Morrey spaces and weighted $L^2$ spaces
with power weights \cite{S3,S4}.

Recently in \cite{S5}
we developed a new approach to the $L^p$ estimates
for boundary value problems,
via $L^2$ estimates, reverse H\"older inequalities,
 and a real variable argument.
For the second order elliptic systems as well as the polyharmonic
equation, we were able to show that the $L^p$ Dirichlet
problem is uniquely solvable for 
$$
n\ge 4\ \ \ \text{ and }\ \ \ \ 
2-\e<p<\frac{2(n-1)}{n-3}
+\e.
\tag 1.7
$$
We remark that in the case of the polyharmonic equation $\Delta^\ell u=0$,
the range in (1.7) is known to be
sharp for $\ell\ge 2$ and $4\le n\le 2\ell +1$.
This was pointed out by Pipher and Verchota 
 \cite{PV1,PV3,PV4}, using examples in \cite{MNP,KM}.

In this paper we continue the work in \cite{S5}
and study the general higher order elliptic
equations and systems. Let
$\Delta(Q,r)=B(Q,r)\cap \partial\OO$ and
$T(Q,r)=B(Q,r)\cap \OO$ where $Q\in \partial\OO$
and $0<r<r_0$.
One of the key ingredients in the approach
we developed in \cite{S5} is the following
(weak) reverse H\"older inequality
$$
\left(\frac{1}{r^{n-1}}\int_{\Delta(Q,r)}
|(\nabla ^{\ell-1} v)^*|^p\, d\sigma\right)^{1/p}
\le C\, \left(\frac{1}{r^{n-1}}\int_{\Delta(Q,2r)}
|(\nabla^{\ell-1}v)^*|^2\, d\sigma\right)^{1/2},
\tag 1.8
$$ for solutions of $\CL(D) v=0$ in 
$\OO$ satisfying $(\nabla^{\ell-1} v )^*\in L^2(\partial\OO)$
and $D^\alpha v=0$ for $|\alpha|\le \ell -1$ on $\Delta(Q,3r)$.
We will show in this paper that 
given a general system of elliptic operators
 $\CL(D)$, a bounded Lipschitz domain
$\OO$ and $p>2$, the reserve H\"older condition
(1.8) with exponent $p$ for $L^2$ solutions with 
zero Dirichlet data on $\Delta(Q,3r)$
is  {\it necessary and sufficient} for the
solvability of the $L^p$ Dirichelt problem on $\OO$.

\proclaim{\bf Theorem 1.9}
Let $\CL(D)$ be a system of elliptic operators of order $2\ell$
 given by (1.1) and satisfying
conditions (1.2) and (1.3). For any bounded Lipschitz domain $\OO$ and $p>2$,
the following are equivalent.

\item{(i)} The $L^p$ Dirichlet problem for $\CL(D) u=0$ on $\OO$ is uniquely
solvable.

\item{(ii)} There exist $C>0$ and $r_0>0$ such that for any $Q\in\partial\OO$
and $0<r<r_0$, the reverse H\"older condition (1.8) holds for
any solution $v$ of $\CL(D) v=0$ in $\OO$ with the properties
$(\nabla^{\ell-1} v)^*\in L^2(\partial\OO)$ and $D^\alpha v=0$ 
for $|\alpha|\le \ell -1$ on $\Delta(Q,3r)$.
\endproclaim

Since the reverse H\"older condition (1.8) has the self-improving
property, it follows that the set of all exponents $p$ in $(2,\infty)$
for which the $L^p$ Dirichlet problem for $\CL(D)u=0$ on $\OO$ is solvable
is an open interval $(2,q)$ with $2<q\le \infty$.

Using square function estimates for $\CL(D)$ \cite{DKPV} as well as 
the regularity estimate (1.12) \cite{PV3,V4}, we also obtain a
much simpler condition which implies (1.8) (see condition (1.11)
as well as (2.15)).
This leads to the following result.

\proclaim{Theorem 1.10}
Let $\CL(D)$ be a system of elliptic operators given by (1.1) and 
 satisfying the symmetry condition (1.2) and
 ellipticity condition (1.3). Let $\OO$ be a bounded
Lipschitz domain in $\br^n$, $n\ge 4$.
Suppose that there exist constants $C_0>0$, $R_0>0$, 
and $\lambda\in (0,n]$ such that
for $0<r<R< R_0$ and $Q\in\partial\OO$,
$$
\int_{T(Q,r)} |\nabla^{\ell-1} v|^2\, dx
\le C_0\, \left(\frac{r}{R}\right)^\lambda \int_{T(Q,R)}
|\nabla^{\ell-1} v|^2\, dx ,
\tag 1.11
$$
whenever
 $v$ is a solution of $\CL(D){v} =0$ in $\OO$ with the properties,
$(\nabla^{\ell -1} v)^*\in L^2(\partial\OO)$ and
$D^\alpha v=0$ on $\Delta(Q, R)$ for 
 $ |\alpha|\le \ell -1$.
Then, if
$$
2<p<2+\frac{4}{n-\lambda},
\tag 1.12
$$
the $L^p$ Dirichlet problem (1.4) on $\OO$ is uniquely solvable.
\endproclaim

 For solutions of the higher order elliptic
equations and systems $\CL(D) u=0$ in $\OO$,
 the following regularity estimate,
$$
\| (\nabla^{\ell} u)^*\|_{L^p(\partial\OO)}
\le C\, \| \nabla_t \nabla^{\ell-1} u\|_{L^p(\partial\OO)},
\tag 1.13
$$
was established by Pipher and Verchota 
 \cite{PV3,V4} for $n\ge 2$ and $2-\e_1<p<2+\e_1$. 
Using (1.13), it is not hard to show that condition (1.11)
holds for some $\lambda >3$.
 As a consequence, we obtain
 the following extension of the main results in \cite{S5}. 
It gives an affirmative answer to an open question raised in
\cite{PV3}.

\proclaim{Corollary 1.14}
For a general higher order homogeneous elliptic system with
real constant coefficients
 satisfying the
symmetry condition and the Legendre-Hadamard ellipticity condition,
the $L^p$ Dirichlet problem (1.4) on $\OO$ is uniquely solvable
for $n\ge 4$ and $2-\e< p<\frac{2(n-1)}{n-3} +\e$, where
$\e>0$ depends on $n$, $m$, $\ell$, $\mu$ and the
Lipschitz character of $\OO$.
\endproclaim

Whether condition (1.11) is necessary for the conclusion of Theorem 1.10
remains open for $p>2(n-1)/(n-3)$ (see Remark 5.21).
As we mentioned earlier,
for the polyharmonic equation $\Delta^\ell u=0$ in $\OO$ 
where $\ell\ge 2$, the $L^p$ Dirichlet problem
is not solvable in general
on Lipschitz domains if $4\le n\le 2\ell +1$ and $p>2(n-1)/(n-3)$.
Thus the range of $p$ in Corollary 1.13 is sharp
in the case $4\le n\le 2\ell +1$.
If $n\ge 2\ell +2$, the $L^p$ Dirichlet problem
is known to be not solvable in general
for $p>2\ell/(\ell -1)$ \cite{PV3,MNP}.
Note that if $\CL(D)$ on $\OO$ satisfies condition (1.11) for some
$\lambda>n-2\ell+2$, which would imply that
$v$ is H\"older continuous up to the boundary $\Delta(Q,R)$,
then the $L^p$ Dirichlet problem is indeed solvable
for $2-\e<p<2\ell/(\ell-1) +\e$
by Theorem 1.10.
In this regards, the problem seems to be closely related to 
the Wiener's type characterization
of regularity for higher order elliptic equations
studied in \cite{M1,MN,M2,M3}.
It is not hard to see
 that for the subclass of the operators
$\CL(D)$ called positive with the weight $F$ studied by
Maz'ya in \cite{M2,M3}
($F$ is the fundamental solution of $\CL (D)$), estimate
(1.11) holds for some $\lambda>n-2\ell+2$ on Lipschitz domains.
See Theorem 2.14.
In particular,
 estimate (1.11) for some $\lambda>n-2\ell+2$ holds
if $\ell=2$ (the biharmonic equation) and $n=5$, $6$, or $7$
\cite{M1}.
In the case $\ell\ge 3$, estimate (1.11) holds for
some $\lambda>n-2\ell+2$ if $n=2\ell +1$ or $2\ell +2$ \cite{MN}.   
This, combined with Theorem 1.10 as well as
results in \cite{DKV1,V3,PV1,PV2,PV4,S5}, yields the
following.

\proclaim{Theorem 1.15}
For the biharmonic equation $\Delta^2 u=0$ in $\OO$, the 
$L^p$ Dirichlet problem is uniquely solvable if
$$
\alignedat3
&2-\e<p\le \infty & \quad\quad\quad &\text{ if }\ \  n=2,\ 3,\\
&2-\e <p< 6+\e & \quad\quad\quad &\text{ if }\ \  n=4,\\
&2-\e <p< 4+\e & \quad\quad\quad &\text{ if }\ \  n=5,\ 6, \ 7,\\
&2-\e< p< \frac{2(n-1)}{n-3} +\e & \quad\quad\quad &\text{ if }\ \  n\ge 8.
\endalignedat
$$
The ranges of $p$ are sharp for $2\le n\le 7$.
\endproclaim

\proclaim{Theorem 1.16}
For the polyharmonic equation $\Delta^\ell u=0$ in $\OO$ with
$\ell \ge 3$, the $L^p$ Dirichlet problem is uniquely
solvable if
$$
\alignedat3
 & 2-\e <p\le \infty & \quad\quad\quad &\text{ if }\ \  n=2,\ 3,\\
& 2-\e <p <\frac{2(n-1)}{n-3} +\e & \quad\quad\quad &\text{ if }\ \ 
4\le n\le 2\ell +1 \text{ or } \ n\ge 2\ell +3,\\
& 2-\e <p< \frac{2\ell}{\ell -1} +\e &\quad\quad\quad & \text{ if }
\ \  n=2\ell +2.
\endalignedat
$$
The ranges of $p$ are sharp for $2\le n\le 2\ell +2$.
\endproclaim

We should remark that only the case $n=6,7$ in Theorem 1.15 and the case
$n=2\ell +2$ in Theorem 1.16 are new.
Also the problem of sharp ranges of $p$ remains open
for $n\ge 8$ if $\ell =2$, and
for $n\ge 2\ell +3$ if $\ell\ge 3$.

The paper is organized as follows.
In Section 2 we collect some basic estimates for solutions of
the elliptic system $\CL(D)u=0$ as well as some inequalities
that will be used in later sections.
 In Section 3 we show that the reverse H\"older condition (1.8)
is sufficient for the $L^p$ solvability of the Dirichlet problem (1.4)
(see Theorem 3.2).
That this condition  is also necessary
is proved in Section 4 (see Theorem 4.1).
While Theorem 1.9 follows by combining Theorem 3.2 with Theorem 4.1,
the proof of Theorem 1.10 as well as Corollary 1.14
will be given in Section 5.
Finally we remark that we will make no effort to distinguish
 vector-valued functions from real valued functions.
It should be clear from the context.

\medskip

\centerline{\bf 2. Basic Estimates and Inequalities}

Most of the materials in this section are known.

\proclaim{\bf Proposition 2.1} (Interior estimates)
Let $u$ be a solution of $\CL(D)u=0$ in $\OO$. Suppose
$B(x,2r)\subset\OO$. Then
$$
|D^\alpha u(x)|\le \frac{C_\alpha}{r^{n+|\alpha|}}
\int_{B(x,r)} |u(y)|\, dy,
\tag 2.2
$$
for any multi-index $\alpha$.
\endproclaim

Estimate (2.2) is well known (see e.g. \cite{DN}). 
In the case when $n$ is odd,
it follows from the potential representation of $u\varphi$ by the fundamental
solution homogeneous of degree $2\ell m-n$
for the elliptic operator det$(\CL^{jk}(D))_{m\times m}$.
If $n$ is even and $2\ell m\ge n$,
 the fundamental solution contains the logarithmic
function $\ln |x|$ (see e.g. \cite{H}, p.169).
 However in this case, one may reduce the problem to the
odd dimensional case by adding an independent
variable (the method of descending).

Let $\psi:\br^{n-1}\to\br$ be a Lipschitz function such that
$\psi(0)=0$. For $r>0$, define
$$
\align
I(r)& =\big\{ (x^\prime,\psi(x^\prime))\in \br^{n-1}\times\br:\
|x_1|<r,\dots, |x_{n-1}|<r\big\}, \tag 2.3\\
Z(r)& =\big\{ (x^\prime, x_n): \
|x_1|<r,\dots, |x_{n-1}|<r,\, \psi(x^\prime)
<x_n<C_0\,  r\big\},\tag 2.4
\endalign
$$
where $C_0=1+10\sqrt{n}\|\nabla \psi\|_\infty>0$ is chosen so that
 $Z(r)$ is a star-shaped
Lipschitz domain with Lipschitz constant independent of $r$.

\proclaim{\bf Lemma 2.5} (Poincar\'e inequality)
Let $1\le p<\infty$. Suppose that $u\in W^{1,p}(Z(r))$ and
$u=0$ on $I(r)$. Then
$$
\int_{Z(r)}|u|^p\, dx
\le C\, r^p\int_{Z(r)} |\nabla u|^p\, dx,
\tag 2.6
$$
where $C$ depends only on $\|\nabla\psi\|_\infty$, $p$ and $n$.
\endproclaim

\demo{Proof} The case $p=1$ follows easily from the fundamental
theorem of calculus. For $p>1$, since $|\nabla |u|^p|
\le p |u|^{p-1}
|\nabla u|$, one applies inequality (2.6)
with $p=1$ to the function $|u|^p$ and then use the H\"older's inequality.
\enddemo

\proclaim{\bf Lemma 2.7} (Sobolev inequality)
Let $1\le p<n$ and $q=pn/(n-p)$. Suppose $u\in W^{1,p}(Z(r))$ and
$u=0$ on $I(r)$. Then
$$
\left(\int_{Z(r)} |u|^q\, dx\right)^{1/q}
\le C\, \left(\int_{Z(r)} |\nabla u|^p\, dx\right)^{1/p},
\tag 2.8
$$
where $C$ depends only on $\|\nabla\psi\|_\infty$, $p$ and $n$.
\endproclaim
\demo{Proof}
Estimate (2.8) follows from (2.6) and the well known Sobolev inequality
$$
\| u\|_{L^q(Z(r))}
\le C\, \bigg\{ \|\nabla u\|_{L^p(Z(r))}
+\frac{1}{r}\, \| u\|_{L^p(Z(r))}\bigg\}.
$$
\enddemo

\proclaim{\bf Lemma 2.9} (Cacciopoli's inequality)
Let $\OO=Z(3r)$. Suppose $\CL(D) u=0$ in $\OO$ and $(\nabla^\ell
 u)^*\in L^2(\partial\OO)$.
Also assume that $D^\alpha u=0$ on $I(3r)$ for $|\alpha|\le \ell-1$.
Then
$$
\int_{Z(r)} |\nabla^\ell u|^2\, dx
\le\frac{C}{r^2}
\int_{Z(2r)}
|\nabla ^{\ell-1} u|^2\, dx,
\tag 2.10
$$
where $C$ depends only on $\|\nabla\psi\|_\infty$, the ellipticity
constant $\mu$ as well as $n$, $m$ and $\ell$.
\endproclaim

\demo{Proof} Let $\varphi$ be a smooth cut-off function on $\br^n$
such that $\varphi=1$ on $Z(r)$, $\varphi=0$ on $Z(3r)\setminus Z(2r)$
and $|D^\alpha \varphi|\le C/r^{|\alpha|}$ for $|\alpha|\le 2\ell$.
To prove (2.10), we use the test function $u\varphi^2$ and
proceed as in the proof of the Cacciopoli's inequality for second order
elliptic systems (\cite{Gr}, pp.76-77).
This gives
$$
\int_{Z(r)}
|\nabla^\ell u|^2\, dx
\le C\, \sum_{|\alpha|\le \ell-1}
\frac{1}{r^{2\ell-2|\alpha|}}
\int_{Z(2r)}
|D^\alpha u|^2\, dx.
$$
We remark that with the assumption $(\nabla^{\ell} u)^*
\in L^2(\partial\OO)$, the necessary integration by parts may be justified
by approximating $\OO$ from inside by a sequence of smooth domains \cite{V1}.
The desired estimate (2.10) now follows from Poincar\'e inequality
(2.6).
\enddemo

\proclaim{\bf Lemma 2.11} (Higher integrability)
Under the same assumption as in Lemma 2.9, we have
$$
\left(\frac{1}{r^n}\int_{Z(r)} |\nabla^{\ell} u|^q\, dx\right)^{1/q}
\le C\, \left(\frac{1}{r^n}\int_{Z(2r)}
|\nabla^\ell u|^2\, dx\right)^{1/2},
\tag 2.12
$$
where $q>2$ depends only on $\|\nabla \psi\|_\infty$, $\mu$, $n$, $m$
and $\ell$.
\endproclaim

\demo{Proof} It follows from Cacciopoli's inequality (2.10)
and Sobolev inequality (2.8) that
$$
\left(\frac{1}{r^n}\int_{Z(r)}
|\nabla^{\ell} u|^2\, dx\right)^{1/2}
\le C\, \left(\frac{1}{r^n}\int_{Z(2r)}
|\nabla^\ell u|^{p_n}\, dx\right)^{1/p},
\tag 2.13
$$
where $p_n=2n/(n+2)$. It is well known that boundary reverse H\"older
inequality (2.13), together with the interior estimate (2.2),
implies  the inequality (2.12) (see \cite{Gr}, pp.122-123).
\enddemo

We end this section with a theorem concerning the condition (1.11).
Recall that for $Q\in \partial\OO$ and $r>0$,
$\Delta(Q,r)= B(Q,r)\cap\partial\OO$
and $T(Q,r)=B(Q,r)\cap \OO$.

\proclaim{\bf Theorem 2.14}
Let $\OO$ be a bounded Lipschitz domain in $\br^n$.
Fix $Q\in \partial\OO$ and $R_0>0$ sufficiently small.
 Let $u$ be a solution of
$\CL(D) u=0$ in $\OO$ with the properties $(\nabla^\ell u)^*
\in L^2(\Delta(Q,R_0))$ and $D^\alpha u=0$ on $\Delta(Q,R_0)$
 for $|\alpha|\le \ell-1$.
Suppose that for some $\lambda_0>0$ and all $0<r<R<R_0$,
$$
\int_{T(Q,r)}
|u|^2\, dx
\le C \, \left(\frac{r}{R}\right)^{\lambda_0 +2\ell-2}
\int_{T(Q,R)} |u|^2\, dx.
\tag 2.15
$$
Then if $0<\lambda<\lambda_0$, we have
$$
\int_{T(Q,r)}|\nabla^{\ell-1} u|^2\, dx
\le C_\lambda\, \left(\frac{r}{R}\right)^\lambda \int_{T(Q,R)}
|\nabla^{\ell-1} u|^2\, dx,
\tag 2.16
$$
for all $0<r<R<R_0$.
\endproclaim

\demo{Proof} By translation and rotation,
 it suffices to prove the theorem with $\Delta(Q,R_0)$,
$T(Q,r)$ and $T(Q,R)$ replaced by
$I(R_0)$, $Z(r)$ and $Z(R)$ respectively.
We may also assume that $0<r<R/2<R_0/4$.

By the interpolation inequality (\cite{A}, p.79)
and Poincar\'e inequality (2.6), we have
$$
\aligned
\| \nabla^{\ell-1} u\|_{L^2(Z(r))}
& \le C\, \| \nabla^\ell u\|_{L^2(Z(r))}^{\frac{\ell-1}{\ell}}
\| u\|_{L^2(Z(r))}^{\frac{1}{\ell}}\\
&\le C\, \| \nabla^{\ell-1} u\|_{L^2(Z(2r))}^{\frac{\ell-1}{\ell}}
\cdot \frac{1}{r^{\frac{\ell-1}{\ell}}}
\cdot \left(\frac{r}{R}\right)^{\frac{\lambda_0 +2\ell -2}{2\ell}}
\| u\|_{L^2(Z(R))}^{\frac{1}{\ell}},
\endaligned
\tag 2.17
$$
where we also use Cacciopoli's inequality (2.10) and assumption (2.15).
It then follows from Poincar\'e inequality (2.6) 
and H\"older's inequality that
$$
\aligned
\|\nabla^{\ell-1} u\|_{L^2(Z(r))}
&\le C\, \|\nabla^{\ell-1} u\|_{L^2(Z(2r))}^{\frac{\ell-1}{\ell}}
\cdot \left(\frac{r}{R}\right)^{\frac{\lambda_0}{2\ell}}
\| \nabla^{\ell-1} u\|_{L^2(Z(R))}^{\frac{1}{\ell}}\\
&\le \varepsilon\, \| \nabla^{\ell-1} u\|_{L^2(Z(2r))}
+C_\varepsilon\left(\frac{r}{R}\right)^{\frac{\lambda_0}{2}}
\|\nabla^{\ell-1} u\|_{L^2(Z(R))},
\endaligned
\tag 2.18
$$
for any $\varepsilon>0$ and any $0<r<R/2$.
By Lemma 8.23 in \cite{GT} (p.201),  estimate (2.18) implies 
$$
\int_{Z(r)} |\nabla^{\ell-1} u|^2\, dx
\le C_\lambda \, \left(\frac{r}{R}\right)^\lambda
\int_{Z(R)} |\nabla^{\ell-1} u|^2\, dx
\tag 2.19
$$
for any $\lambda\in (0,\lambda_0)$ and for all $0<r<R<R_0$.
The proof is finished.
\enddemo

\medskip

\centerline{\bf 3. The Sufficiency of the Reverse H\"older Condition}

The goal of this section is to show that
for a given operator $\CL(D)$ on a fixed Lipschitz
domain $\OO$, the reverse H\"older condition
(1.8) with exponent $p$
 for solutions with zero Dirichlet data on $\Delta(Q,3r)$ 
is sufficient for the solvability of the $L^p$ Dirichlet problem on $\OO$.

Recall that for a function $u$ defined on $\OO$, 
its non-tangential maximal function
$(u)^*$ is defined by
$$
(u)^*(Q)
=\sup\big\{
|u(x)|:\  x\in\Gamma(Q)\big\}, \ \ \ \ \text{ for } Q\in \partial\OO,
\tag 3.1
$$
where $\Gamma(Q)
=\big\{ x\in \OO: \, |x-Q|<2\, \text{dist}(x,\partial\OO)\big\}$.

\proclaim{\bf Theorem 3.2}
Let $\CL(D)$ be an elliptic operator given by (1.1) and satisfying the 
symmetry condition (1.2) and ellipticity condition (1.3).
Let $\OO$ be a bounded Lipschitz domain in $\br^n$.
Fix $p>2$. Suppose that for any $\Delta(Q,r)\subset\partial\OO$,
 inequality (1.8) holds for
all solutions of $\CL(D)v=0$ in $\OO$ with the properties
$(\nabla^{\ell-1} v)^*\in L^2(\partial\OO)$ and
$D^\alpha v=0$ for $|\alpha|\le \ell-1$ on $\Delta(Q,3r)$.
Then the $L^p$ Dirichlet problem (1.4) is uniquely solvable.
\endproclaim

Since $p>2$, the uniqueness in Theorem 3.2 follows from the
uniqueness for the case $p=2$ \cite{V4}. For the existence as well as
estimate (1.6), it suffices to establish the following lemma.

\proclaim{\bf Lemma 3.3}
Let $(\dot{f},g)\in \big(WA^p_{\ell -1}(\partial\OO), L^p(\partial\OO)\big)$.
Let $u$ be the unique $L^2$ solution of (1.4) with boundary data
$(\dot{f},g)$, i.e., $u$ satisfies
(1.4) and $(\nabla^{\ell -1} u)^* \in L^2(\partial\OO)$.
Then, under the same assumption as in Theorem 3.2, we have
$$
\aligned
&\left\{\frac{1}{r^{n-1}}
\int_{\Delta(Q,r)}
|(\nabla^{\ell-1} u)^*|^p\, d\sigma\right\}^{1/p}
\le C\, \left\{\frac{1}{r^{n-1}}
\int_{\Delta(Q,2r)}
|(\nabla^{\ell-1} u)^*|^2\, d\sigma\right\}^{1/2}\\
&\ \ \ \ \ \ \ \ \ \ \ \ \ \ \ \ \ \ \ \ \ \ 
+C\, \bigg\{\frac{1}{r^{n-1}}
\int_{\Delta(Q,2r)}
\big(|g|^p + \sum_{|\alpha|=\ell-2}
|\nabla_t f_\alpha|^p\big)\, d\sigma\bigg\}^{1/p},
\endaligned
\tag 3.4
$$
for any $\Delta(Q,r)\subset\partial\OO$ with $Q\in \partial\OO$ and $0<r<r_0$.
\endproclaim

The desired estimate (1.6) follows from (3.4)
by covering $\partial\OO$ with a finite number of 
surface balls $\Delta(Q,r)$.

The proof of Lemma 3.3 is essentially contained in \cite{S5}, where
the cases of second order elliptic systems and the polyharmonic
equation were considered. For completeness as well as reader's
convenience we will provide
a detailed proof here. 

We will need the following
 Poincar\'e type inequality on $\partial\OO$.

\proclaim{\bf Lemma 3.5}
Suppose $\ell\ge 2$. Let $\dot{f}=\{ f_\alpha:\, |\alpha|\le\ell -2\}
\in WA^2_{\ell -1}(\partial\OO)$
and $\Delta(Q,r)\subset \partial\OO$. Then there exists a polynomial
$h$ of degree at most $\ell-2$ such that
$$
\| f_\beta -D^\beta h \|_{L^2(\Delta(Q,r))}
\le C\, r^{\ell -|\beta|-1}
\sum_{|\alpha|=\ell -2}
\|\nabla_t f_\alpha\|_{L^2(\Delta(Q,r))},
\tag 3.6
$$
for any multi-index $\beta$ with $|\beta|\le \ell -2$.
\endproclaim

\demo{Proof}
Let $h(x)=\sum_{|\alpha|\le \ell -2} \frac{C_\alpha}{\alpha !} \, x^\alpha$ be
a  polynomial where 
 $C_\alpha$ is defined inductively by
$$
\aligned
&C_\alpha =\frac{1}{|\Delta(Q,r)|}
\int_{\Delta(Q,r)} f_\alpha\, d\sigma \ \ \ \text{ if } |\alpha|=\ell -2,\\
&C_{\alpha}
=\frac{1}{|\Delta(Q,r)|}
\int_{\Delta(Q,r)}
\bigg\{ f_\alpha (P)-\sum\Sb  \beta>\alpha\\ |\beta|\le \ell -2\endSb
\frac{C_\beta}{(\beta-\alpha)!}\, P^{\beta-\alpha}\bigg\}\, d\sigma(P),
\ \ \text{ if } |\alpha|<\ell-2.
\endaligned
$$
It is easy to check that
$$
\int_{\Delta(Q,r)}
\big\{ f_\alpha -D^\alpha h\big\}\, d\sigma
=0\ \ \ \text{ for all } |\alpha|\le \ell -2.
\tag 3.7
$$
With (3.7), by using Poincar\'e inequality on $\Delta(Q,r)$ repeatedly, 
we obtain
$$
\aligned
\| f_\beta -D^\beta h\|_{L^2(\Delta(Q,r))}
& \le C\, r\, \sum_{|\alpha|=|\beta|+1}
\| f_\alpha -D^\alpha h\|_{L^2(\Delta(Q,r))}\\
&\le C\, 
r^{\ell -|\beta|-2}
\sum_{|\alpha|=\ell-2}
\| f_\alpha -D^\alpha h\|_{L^2(\Delta(Q,r))}\\
&\le C\, 
r^{\ell -|\beta|-1}
\sum_{|\alpha|=\ell-2}
\|\nabla_t f_\alpha\|_{L^2(\Delta(Q,r))},
\endaligned
$$
for any $|\beta|\le \ell -2$. The proof is finished.
\enddemo

To prove (3.4), we fix
$\Delta(Q,r)$ with $Q\in \partial\OO$ and $0<r<r_0$.
By rotation and translation, we may assume that
$Q=0$ and 
$$
\aligned
B(0,C_1\,r_0)\cap\OO
&=B(0,C_1 \, r_0)\cap \big\{(x^\prime, x_n)\in \br^{n-1}\times\br:\
x_n>\psi(x^\prime)\big\},\\
B(0,C_1\, r_0)\cap\partial\OO
&=B(0, C_1\, r_0)\cap
\big\{ (x^\prime, \psi(x^\prime)):\ x^\prime\in \br^{n-1}\big\},
\endaligned
\tag 3.8
$$
where $\psi$ is a Lipschitz function on $\br^{n-1}$.
Let $S=\big\{ (x^\prime, \psi(x^\prime)):\ x^\prime \in \br^{n-1}\big\}$.
We will perform a Calder\'on-Zygmund decomposition on $S$ in the
proof of (3.4). To do this, we need to introduce
surface cubes on the set $S$.

Let $\Phi: S\to \br^{n-1}$ be a map defined by $\Phi(x^\prime,\psi(x^\prime))
=x^\prime$. A subset $I$ of $S$ is said to be a cube on $S$
if $\Phi(I)$ is a cube in $\br^{n-1}$ with sides parallel
to the coordinate planes. A cube $I$ on $S$ is said to be a
dyadic subcube of $I^\prime$ if $\Phi(I)$ is a dyadic subcube
of $\Phi(I^\prime)$.
Also for $\rho>0$ and a cube $I$ on $S$, we will use $\rho I$ to denote
$\Phi^{-1}(\rho\Phi(I))$.

For cube $I$ on $S$ and a function $f$ defined on $I$, we
define a localized Hardy-Littlewood maximal function
$M_I$ by
$$
M_I(f)(P)
=\sup\Sb \text{cube }I^\prime\owns P\\ I^\prime\subset I\endSb
\frac{1}{|I^\prime|}
\int_{I^\prime}
|f|\, d\sigma\ \ \ \ \ \text{ for } P\in I.
\tag 3.9
$$ 

\demo{\bf Proof of Lemma 3.3} The proof relies on a real
variable argument which is motivated by the method used in \cite{CP}.
Let $I$ be a cube on $S$ such that
$3I\subset S\cap B(0,2r_0)$.
For $\lambda>0$, let
$$
E(\lambda)=\big\{ Q\in I:\
M_{2I}(|(\nabla^{\ell -1} u)^*|^2)(Q)>\lambda\big\}.
\tag 3.10
$$
Fix $2<q<p$, let $A=1/(2\delta^{2/q})$ where $\delta>0$ is a small constant
to be determined. 
Let $F=|g|^2 +\sum_{|\alpha|=\ell-2} |\nabla_t f_\alpha|^2$.
We claim that there exist positive constants
$\delta$, $\gamma$ and $C_0$ depending only on
$n$, $m$, $\mu$, $\ell$, $\OO$ as well as the constant $C$ in the
reverse H\"older condition (1.8) such that
$$
|E(A\lambda)|
\le \delta |E(\lambda)|
+|\big\{
Q\in I:\ M_{2I} 
(F)(Q)>\gamma \lambda\big\}|
\tag 3.11
$$
for all $\lambda\ge \lambda_0$, where
$$
\lambda_0 =\frac{C_0}{|2I|}
\int_{2I}|(\nabla^{\ell -1} u)^*|^2\, d\sigma.
\tag 3.12
$$

Assume the claim is true for a moment. We multiply both sides
of (3.11) by $\lambda^{\frac{q}{2}-1}$ and integrate the
resulting inequality in $\lambda$ over the interval $(\lambda_0, \Lambda)$.
We obtain
$$
\frac{1}{A^{\frac{q}{2}}}
\int_{A\lambda_0}^{A\Lambda}
\lambda^{\frac{q}{2}-1}
|E(\lambda)|\, d\lambda
\le \delta\int_{\lambda_0}^{\Lambda}
\lambda^{\frac{q}{2}-1}|E(\lambda)|\, d\lambda
+C\int_{2I}
|F|^{\frac{q}{2}}\, d\sigma,
\tag 3.13
$$
where we have used the boundedness of $M_{2I}$ on $L^{q/2}(2I)$.
This implies that
$$
\aligned
\left(\frac{1}{A^{q/2}} -\delta\right)
&\int_0^\Lambda \lambda^{\frac{q}{2}-1}
|E(\lambda)|\, d\lambda\\
&\le C\, 
\int_0^{A\lambda_0}
\lambda^{\frac{q}{2}-1} |E(\lambda)|\, d\lambda
+C\, 
\int_{2I}
|F|^{\frac{q}{2}}\, d\sigma\\
&\le C\lambda_0^{\frac{q}{2}}
|I|
+C\, \int_{2I}
\bigg\{
|g|^q +\sum_{|\alpha|=\ell-2}|\nabla_t f_\alpha|^q\bigg\}\, d\sigma.
\endaligned
$$
Observe that $\delta A^{q/2}<1$. Let $\Lambda\to \infty$ in the above
inequality, we obtain $(\nabla^{\ell-1} u)^*\in L^q(I)$ and
$$
\aligned
\bigg\{\frac{1}{|I|}
\int_I |(\nabla^{\ell-1} u)^*|^q\, d\sigma\bigg\}^{1/q}
&\le C \bigg\{\frac{1}{|2I|}
\int_{2I} |(\nabla^{\ell-1} u)^*|^2\, d\sigma\bigg\}^{1/2}\\
&\ \ \ \ \ \
+C\bigg\{ \frac{1}{|2I|}\int_{2I}
\big(|g|^q +\sum_{|\alpha|=\ell-2}
|\nabla_t f_\alpha |^q \big)\, d\sigma\bigg\}^{1/q},
\endaligned
\tag 3.14
$$
for any $q\in (2,p)$.  
Note that the reverse H\"older condition (1.8)
is a self-improving property. That is, if  $(\nabla^{\ell-1}u)^*$
satisfies condition (1.8) for some $p>2$, then it
satisfies condition (1.8) for some $\bar{p}>p$ (\cite{Gi}, pp.122-123).
Thus (3.14) holds for $q=p$.
Estimate (3.4) now follows from (3.14) by covering $\Delta(Q,r)$
with a finite number of sufficiently small surface cubes $I$. 

It remains to prove the claim (3.11). To this end, we fix
$\lambda\ge \lambda_0$. Note that
$E(\lambda)$ is open relative to $I$. This implies that there exists
a sequence of disjoint dyadic subcubes $\{ I_j\}$ of $I$ such that
$E(\lambda)=\cup I_j$ (up to a set of surface measure zero).
We may assume that each $I_j$ is maximal in the sense that
if $I^\prime \supset I_j$ is a dyadic subcube of $I$, then
$I^\prime\nsubseteq E(\lambda)$ unless $I^\prime=I_j$.
Also, since
$$
|E(\lambda)|\le \frac{C}{\lambda}
\int_{2I} |(\nabla^{\ell-1} u)^*|^2\, d\sigma
\le\frac{C\lambda_0 |I|}{C_0\lambda}
\le \frac{C |I|}{C_0},
\tag 3.15
$$
we may assume that $|E(\lambda)|\le \delta |I|$ by taking
$C_0$ sufficiently large. It follows that $|I_j|\le \delta |I|$.
In particular we may assume that $32I_j\subset 2I$.

We claim that it is possible to choose positive constants
$\delta$, $\gamma$ and $C_0$ such that if
$$
\big\{ Q\in I_j:\ M_{2I}(F)(Q)\le \gamma\lambda\big\}
\neq\emptyset,
\tag 3.16
$$
then $|E(A\lambda)\cap I_j|\le \delta |I_j|$.
Clearly, this yields estimate (3.11) by summation.

To prove the last claim, we fix $I_j$ which satisfies (3.16).
Note that for any $Q\in I_j$,
$$
M_{2I} (|(\nabla^{\ell-1} u)^*|^2)(Q)
\le \max \big(M_{2I_j}(|(\nabla^{\ell-1} u)^*|^2)(Q), C_1\lambda\big)
$$
where $C_1$ depends only on $n$ and $\|\nabla\psi\|_\infty$.
This is because $I_j$ is maximal. We may assume that $A\ge C_1$
by taking $\delta$ small. It follows that
$$
| E(A\lambda)\cap I_j |=
|\big\{ Q\in I_j: M_{2I_j} (|(\nabla^{\ell-1} u)^*|^2)(Q)>A\lambda\big\}|.
\tag 3.17
$$
Let $\varphi$ be a smooth cut-off function on $\br^n$ such that
$\varphi=1 $ on $16I_j$, $\varphi=0$ on $\partial\OO\setminus
17 I_j$, and $|D^\alpha\varphi|\le C/\rho^{|\alpha|}$
for $|\alpha|\le \ell -1$, where $\rho=\rho_j$ is the diameter of 
$I_j$.
Let $h$ be a polynomial of degree at most $\ell-2$
given by Lemma 3.5, but with $\Delta(Q,r)$ replaced by $17I_j$.
Let $w=w_j$ be the solution of the $L^2$ Dirichlet problem
(1.4) with boundary data
$$
D^\alpha w =D^\alpha ((u-h)\varphi)
=\sum_{\beta\le \alpha}
\frac{\alpha!}{\beta!(\alpha-\beta)!}
(f_\beta -D^\beta h)D^{\alpha-\beta} \varphi
\tag 3.18
$$
for $|\alpha|\le \ell-2$ and $\frac{\partial^{\ell-1} w}{\partial N^{\ell-1}}
=g\varphi$ on $\partial\OO$.
Note that
$$
\aligned
&\sum_{|\alpha|=\ell-2}
\int_{\partial\OO} |\nabla_t D^\alpha w|^2\, d\sigma\\
&\le C\, \bigg\{
\sum_{|\alpha|=\ell-2}
\int_{17I_j} |\nabla_t f_\alpha|^2\, d\sigma
+\sum_{|\beta|\le \ell -2}
\frac{1}{\rho^{\ell-|\beta|-1}}
\int_{17I_j} |f_\beta -D^\beta h|^2\, d\sigma\bigg\}\\
&\le 
C\sum_{|\alpha|=\ell-2}
\int_{17I_j} |\nabla_t f_\alpha|^2\, d\sigma,
\endaligned
\tag 3.19
$$
where we have used Lemma 3.5 in the last inequality. It follows
from the $L^2$ estimates \cite{PV3,V4} and (3.16) that
$$
\aligned
\int_{\partial\OO}
|(\nabla^{\ell-1} w)^*|^2\, d\sigma
&\le C\, \int_{17I_j}
\bigg\{ |g|^2 +\sum_{|\alpha|=\ell-2} |\nabla_t f_\alpha |^2\bigg\}
\, d\sigma\\
&\le C\, \gamma \lambda |I_j|.
\endaligned
\tag 3.20
$$

Now let $v=u-w-h$ in $\OO$. Note that $v$ is a solution to
the $L^2$ Dirichlet problem ( 1.4) with boundary data vanishing
on $16 I_j$. Indeed, $D^\alpha v
=D^\alpha\big((u-h)(1-\varphi)\big)$ 
for $|\alpha|\le \ell-2$ and
$\frac{\partial^{\ell-1} v}{\partial N^{\ell-1}}
=g(1-\varphi)$ on $\partial\OO$.
Hence we may apply the reverse H\"older condition (1.8)
to $v$ on $16I_j$. This gives 
$$
\aligned
&\left(\frac{1}{|2I_j|}\int_{2I_j}
|(\nabla^{\ell-1} v)^*|^p\, d\sigma\right)^{1/p}
\le
C\, \left(\frac{1}{|4I_j|}\int_{4I_j}
|(\nabla^{\ell-1} v)^*|^2\, d\sigma\right)^{1/2}\\
&\le C\, \left(\frac{1}{|4I_j|}
\int_{4I_j}|(\nabla^{\ell-1} u)^*|^2\, d\sigma\right)^{1/2}
+C
\left(\frac{1}{|4I_j|}
\int_{4I_j}|(\nabla^{\ell-1} w)^*|^2\, d\sigma\right)^{1/2}\\
&\le C\,\lambda^{1/2},
\endaligned
\tag 3.21
$$
where we have used the fact
 $\nabla^{\ell-1} u=\nabla^{\ell-1}w +\nabla^{\ell-1}v$ in
$\OO$. In (3.21) we
also use (3.20) as well as the observation $3I_j\nsubseteq E(\lambda)$
for the last inequality.
We should point out that the reverse H\"older condition
on surface balls $\Delta(Q,r)=B(Q,r)\cap\partial\OO$ is equivalent
to the reverse H\"older condition over surface cubes. 
This is because one may cover any surface cube by
sufficiently small surface cubes with a finite overlap
and vice versa.

Finally in view of (3.17), (3.20) and (3.21) we have
$$
\aligned
|E(A\lambda)\cap I_j|
&= |\big\{ Q\in I_j:\
M_{2I_j} (|(\nabla^{\ell-1} u)^*|^2)(Q)>A\lambda\big\}|\\
&\le |\big\{ Q\in I_j:\
M_{2I_j} (|(\nabla^{\ell-1} w)^*|^2)(Q)>\frac{A\lambda}{4}\big\}|\\
&\ \ \ \ \ \ \ \ \ \ \ \ \ \ \
+|\big\{ Q\in I_j:\
M_{2I_j} (|(\nabla^{\ell-1} v)^*|^2)(Q)>\frac{A\lambda}{4}\big\}|\\
&\le
\frac{C}{A\lambda}
\int_{2I_j}
|(\nabla^{\ell-1} w)^*|^2\, d\sigma
+\frac{C}{(A\lambda)^{p/2}}
\int_{2I_j}
|(\nabla^{\ell-1} v)^*|^p\, d\sigma\\
&\le
|I_j|\bigg\{
\frac{C\gamma}{A} +\frac{C}{A^{p/2}}\bigg\}\\
&\le \delta\,|I_j|\big\{
C\delta^{\frac{p}{q}-1}+
C\gamma \delta^{\frac{2}{q}-1}\big\},
\endaligned
$$
where we have used $A=1/(2\delta^{2/q})$ in the last inequality.
Since $q<p$, we may choose $\delta>0$ so small that $C\delta^{\frac{p}{q}-1}
\le 1/2$.
With $\delta$ chosen, we then choose $\gamma>0$ so small that
$C\gamma \delta^{\frac{2}{q}-1}\le 1/2$.
This gives $|E(A\lambda)\cap I_j|\le \delta |I_j|$.
The proof is now complete.
\enddemo

\medskip

\centerline{\bf 4. The Necessity of the Reverse H\"older Condition}

In this section we will show that for a given elliptic operator
$\CL(D)$ on a fixed Lipschitz domain $\OO$,
the reverse H\"older condition
(1.8) with exponent $p$ is also necessary for the
solvability of the $L^p$ Dirichlet problem (1.4).

\proclaim{\bf Theorem 4.1}
Let $\CL(D)$ be an elliptic operator given by (1.1). Let
$\OO$ be a bounded Lipschitz domain in $\br^n$, $n\ge 4$. Fix $p>2$.
Suppose that the $L^p$ Dirichlet problem
(1.4) on $\OO$ is uniquely solvable.
 Then the reverse H\"older inequality
(1.8) holds for solutions of $\CL(D)(v)=0$ in $\OO$ with
the properties $(\nabla^{\ell-1} v)^*\in L^2(\partial\OO)$ and
$D^\alpha v=0$ on $\Delta(Q,3r)$ for $|\alpha|\le \ell-1$.
\endproclaim

To prove Theorem 4.1, we will need a lemma on the traces
of Riesz potentials. Let
$$
I_1 (f)(x)
=\int_\OO \frac{f(y)\, dy}{|x-y|^{n-1}}.
\tag 4.2
$$
\proclaim{\bf Lemma 4.3}
Let $1<q<n$ and $p=q(n-1)/(n- q)$. Then
$$
\| I_1 (f)\|_{L^p(\partial\OO)}
\le C\, \| f\|_{L^q(\OO)}.
\tag 4.4
$$
\endproclaim
\demo{Proof}
Since $\OO$ is a Lipschitz domain, there exists a smooth vector
field $\bold{V}(x)$ on $\br^n$ such that
$\bold{V}\cdot N\ge c_0>0$ on $\partial\OO$ \cite{V1}.
It follows from the divergence theorem that
$$
\aligned
& c_0\int_{\partial\OO} |I_1(f)|^p\, d\sigma
\le \int_{\partial\OO} \bold{V}\cdot N\, |I_1 (f)|^p\, d\sigma\\
&\le C\, \int_{\OO} |I_1 (f)|^p \, dx
+C\, \int_{\OO} |I_1 (f)|^{p-1}|\nabla I_1 (f)|\, dx\\
&\le
C\, \int_{\OO} |I_1 (f)|^p \, dx
+C\, \left(\int_\OO |\nabla I_1 (f)|^{q}\, dx\right)^{1/q}
\left(\int_\OO |I_1 (f)|^{(p-1)q^\prime}\, dx\right)^{1/q^\prime}.
\endaligned
\tag 4.5
$$
Observe that 
$
\frac{1}{(p-1)q^\prime}=\frac{1}{q}-\frac{1}{n}.
$
It follows that the last term on the right side of (4.5)
is bounded by $C\, \| f\|_{L^q(\OO)}^p$.
To see this, one uses the well known estimates for fractional integrals
and singular integrals \cite{St}.
It is clear that the first term on the right side of (4.5) is also
bounded by  $C\, \| f\|_{L^q(\OO)}^p$.
The proof is complete.
\enddemo

\demo{\bf Proof of Theorem 4.1}
Fix $Q_0\in \partial\OO$ and $0<r<r_0$.
Let $v$ be a solution of $\CL(D) v=0$ in $\OO$ with the properties
$(\nabla^{\ell-1}v)^*\in L^2(\partial\OO)$ and $ D^\alpha v=0$ on
$\Delta(Q_0,3r)$ for $|\alpha|\le \ell-1$.
For a function $u$ on $\OO$, define
$$
\aligned
\Cal{M}_1(u)(Q)
&=\sup\big\{ |u(x)|:\ x\in \Gamma(Q)\ \text{ and }\ |x-Q|
<c_0\, r\big\},\\
\Cal{M}_2(u)(Q)
&=\sup\big\{ |u(x)|:\ x\in \Gamma(Q)\ \text{ and }\ |x-Q|
\ge c_0\, r\big\}
\endaligned
\tag 4.6
$$
for $Q\in \partial\OO$. Clearly, $(\nabla^{\ell-1} v)^*=
\max\big\{ \Cal{M}_1(\nabla^{\ell-1} v), \Cal{M}_2(\nabla^{\ell-1} v)\big\}$.
Note that if $x\in\Gamma(Q)$ for some $Q\in\Delta(Q_0,r)$ and 
$|x-Q|\ge c_0\,r$, by interior estimate (2.2), we have
$$
|\nabla^{\ell-1} v(x)|
\le \frac{C}{r^n}\int_{B(x,cr)}
|\nabla^{\ell-1} v(y)|\, dy
\le \frac{C}{r^{n-1}}
\int_{\Delta(Q_0,2r)}
|(\nabla^{\ell-1} v)^*|\, d\sigma.
\tag 4.7
$$
It follows that for any $p>2$,
$$
\left(\frac{1}{r^{n-1}}
\int_{\Delta(Q_0,r)}
|\Cal{M}_2 (\nabla^{\ell-1} v)|^p\, d\sigma\right)^{1/p}
\le C\, \left(\frac{1}{r^{n-1}}
\int_{\Delta(Q_0,2r)}
|(\nabla^{\ell-1} v)^*|^2\, d\sigma\right)^{1/2}.
\tag 4.8
$$

The estimate of $\Cal{M}_1(\nabla^{\ell-1}v)$ on $\Delta(Q_0,r)$
is much more involved.
First, we
choose a smooth cut-off function $\varphi$ on $\br^n$ such that
$\varphi=1$ on $B(Q_0,2r)$, supp$\varphi\subset B(Q_0,3r)$, and
$|D^\alpha \varphi|\le C/r^{|\alpha|}$ for $|\alpha|\le 2\ell$.
Note that 
$$
\aligned
\big[ \CL(D)(v\varphi)\big]^j
&=\sum_{k=1}^m\sum_{|\alpha|=|\beta|=\ell}
a_{\alpha\beta}^{jk}
D^{\alpha+\beta}(v^k\varphi)\\
&=\sum_{k=1}^m\sum_{|\alpha|=|\beta|=\ell}
a_{\alpha\beta}^{jk}
D^\alpha\bigg\{ D^\beta v^k\cdot\varphi
+\sum_{\gamma<\beta}
\frac{\beta!}{\gamma!(\beta-\gamma)!}
D^\gamma v^k\cdot D^{\beta-\gamma} \varphi\bigg\}\\
&=
\sum_{k=1}^m
\sum_{|\alpha|=|\beta|=\ell}
\sum_{\gamma<\alpha}
a_{\alpha\beta}^{jk}
\frac{\alpha!}{\gamma!(\alpha-\gamma)!}
D^{\beta+\gamma}v^k\cdot D^{\alpha-\gamma} \varphi\\
&\ \ \ \ \ \ +\sum_{k=1}^m\sum_{|\alpha|=|\beta|=\ell}
\sum_{\gamma<\beta} a_{\alpha\beta}^{jk}
\frac{\beta!}{\gamma!(\beta-\gamma)!}
D^{\alpha}\big(D^{\gamma}v^k\cdot D^{\beta-\gamma}\varphi\big),
\endaligned
\tag 4.9
$$
where we have used $\CL(D)(v)=0$ in $\OO$.
Let $G(x)=(G_{ij}(x))_{m\times m}$ denote a matrix of fundamental
solutions on $\br^n$ to the operator $\CL(D)$ with pole at the origin.
We remark that if $n$ is odd or $2\ell<n$,
$G_{ij}(x)$ is homogeneous of degree $2\ell-n$ and smooth 
away from the origin.
However, if $n$ is even and $2\ell\ge n$, 
the logarithmic function $\ln |x|$ appears in $G(x)$.
Indeed in this case, we have $G_{ij}(x)
=G_{ij}^{(1)}(x) +G_{ij}^{(2)}(x)\ln |x|$ where
$G_{ij}^{(1)}(x)$ is homogeneous of degree $2\ell-n$ and
$G_{ij}^{(2)}(x)$ is a polynomial of degree $2\ell-n$ (see \cite{F}, p.76
or \cite{H}, p.169). To deal with the factor $\ln |x|$,
 we need to replace $\ln |x|$ by $\ln (|x|/r)$.
This can be done because $G_{ij}^{(2)}(x)$ is a polynomial of
degree $2\ell-n$.

Note that in both cases, we have
$$
|D^\alpha G(x)|\le\frac{C_\alpha}{|x|^{n-2\ell +|\alpha|}}
\ \ \ \ \text{ for } |\alpha|\ge 2\ell-n+1,
\tag 4.10
$$
as  the derivatives $D^\alpha$ 
eliminate the (possible) logarithmic singularity if $|\alpha|> 2\ell-n$.

Next we fix $y_0\in \br^n\setminus\overline{\OO}$ so that
$|y_0-Q_0|=r=\text{dist}(y_0,\partial\OO)$.
Let $\widetilde{G}(x,y)=G(x-y)$ and
$$
F_{ij}(x,y)=\widetilde{G}_{ij}(x,y)
-\sum_{|\gamma|\le 2\ell-1}
\frac{(y-y_0)^\gamma}{\gamma!} D_y^\gamma\widetilde{G}_{ij}(x,y_0).
\tag 4.11
$$
Note that the summation term in (4.11) is a solution to
$\CL(D) u=0$ in $\OO$ in  both $x$ and $y$ variables.
It is added to $\widetilde{G}(x,y)$ in order 
to create the desired decay when $|x-Q_0|\ge 5r$ and
$|y-Q_0|\le 3r$. Indeed,
by the Taylor remainder theorem and (4.10), if
$ x\in \OO\setminus T(Q_0,5r)$ and $ y\in T(Q_0,3r)$, we have
$$
|\nabla_x^{\ell-1} D_y^\alpha F_{ij}(x,y)|
\le \frac{C\, r^{2\ell-|\alpha|}}{|x-y|^{n+\ell-1}}
\ \ \ \ \text{ for } |\alpha|\le 2\ell.
\tag 4.12
$$
Also, if $x\in T(Q_0,5r)$ and $y\in T(Q_0,3r)$,
$$
|\nabla_x^{\ell-1}D_y^\alpha F_{ij}(x,y)|
\le
\frac{C\, r^{\ell-|\alpha|}}
{|x-y|^{n-1}}
\ \ \ \ \text{ for }|\alpha|\le \ell.
\tag 4.13
$$
To see (4.13), one considers two cases:
$|\alpha|>\ell-n+1$ and $|\alpha|\le \ell-n+1$.
In the first case, one uses estimate (4.10).
For the second case, the (possible)
term involving the logarithmic function in
$\nabla_x^{\ell-1} D_y^\alpha \widetilde{G}_{ij} (x,y)$ is bounded by
$C\, |x-y|^{\ell-n-|\alpha|+1} |\ln |\frac{|x-y|}{r}||$.
Since $|x-y|\le C\, r$, it is clearly 
bounded by the right side of (4.13).

 In view of (4.9), we let $w(x)=
(w^1(x),\dots, w^m(x))$ where
$$
\aligned
w^i(x)
&=\sum_{j,k=1}^m\sum_{|\alpha|=|\beta|=\ell}
\sum_{\gamma<\alpha}
(-1)^\ell a_{\alpha\beta}^{jk}
\frac{\alpha!}{\gamma! (\alpha-\gamma)!}
\int_\OO D^\beta_y\big\{
F_{ij}(x,y) D^{\alpha-\gamma} \varphi(y)\big\}
D^\gamma v^k(y)\, dy\\
&\ \ \ \ +
\sum_{j,k=1}^m\sum_{|\alpha|=|\beta|=\ell}
\sum_{\gamma<\beta}
(-1)^\ell a_{\alpha\beta}^{jk}
\frac{\beta!}{\gamma! (\beta-\gamma)!}
\int_\OO D^\alpha_y
F_{ij}(x,y) D^\gamma v^k(y)\cdot D^{\beta-\gamma}
\varphi(y)\, dy.
\endaligned
$$
Then $\CL(D)(w) =\CL(D)(v\varphi)$ in $\OO$. To see this, one may fix
$B(x_0,3s)\subset \OO$ and write $w=w_1+w_2$, where $w_1$ and $w_2$ are
defined in the same way as $w$ but with domain $\OO$
of both integrals replaced by
$B(x_0,2s)$ and $\OO\setminus B(x_0,2s)$ respectively.
 Clearly $\CL(D)(w_2)=0$
in $B(x_0,s)$. To show $\CL(D)(w_1) =
\CL(D)(v\varphi)$ in $B(x_0,s)$, one uses integration by parts and (4.9).

To continue, we observe that on $\Delta(Q_0,r)$, 
$$
\Cal{M}_1(\nabla^{\ell-1} v)
=\Cal{M}_1(\nabla^{\ell-1}(v\varphi))
\le \Cal{M}_1(\nabla^{\ell-1}w)
+\Cal{M}_1(\nabla^{\ell-1}(v\varphi-w)).
\tag 4.14
$$
It follows from (4.13) that for $x\in T(Q_0,5r)$,
$$
|\nabla^{\ell-1} w(x)|
\le C\,
\sum_{|\gamma|\le \ell-1}
\int_{T(Q_0,3r)\setminus T(Q_0,2r)}
\frac{|D^\gamma v(y)|\, r^{|\gamma|-\ell}}
{|x-y|^{n-1}}\, dy.
\tag 4.15
$$
This implies  that if $Q\in \Delta(Q_0,r)$ and
the constant $c_0$ in (4.6) is sufficiently small, 
$$
\aligned
\Cal{M}_1(\nabla^{\ell-1} w)(Q)
&\le C\, \sum_{|\gamma|\le \ell-1}
r^{|\gamma|-\ell-n+1}
\int_{T(Q_0,3r)\setminus T(Q_0,2r)}
|D^\gamma v(y)|\, dy\\
&\le C\, \sum_{|\gamma|\le \ell-1}r^{|\gamma|-\ell+1}
\left(\frac{1}{r^n}\int_{T(Q_0,3r)}
|D^\gamma v(y)|^2\, dy\right)^{1/2}\\
&\le C\, \left(\frac{1}{r^n}\int_{T(Q_0,3r)}
|\nabla^{\ell-1} v|^2\, dy\right)^{1/2},
\endaligned
\tag 4.16
$$
 where we have used
Poincar\'e inequality (2.6) in the last step (the proof for
Poincar\'e inequality on $T(Q,r)$ is the same).
Clearly, this gives
$$
\left(\frac{1}{r^{n-1}}\int_{\Delta(Q_0,r)}
|\Cal{M}_1(\nabla^{\ell-1} w)|^p\, d\sigma\right)^{1/p}
\le C\, \left(\frac{1}{r^{n-1}}
\int_{\Delta(Q_0,3r)}
|(\nabla^{\ell-1} v)^*|^2\, d\sigma\right)^{1/2}.
\tag 4.17
$$

It remains to estimate $\Cal{M}_1(\nabla^{\ell-1}(v\varphi-w))$.
It is here that
 we need to use the assumption that the $L^p$ Dirichlet problem
(1.4) on $\OO$ is uniquely solvable.

Note that $\CL(D)(v\varphi -w)=0$ in $\OO$.
We also have $(\nabla^{\ell-1}(v\varphi -w))^*\in L^2(\partial\OO)$.
To see this, by the square function estimates (see (5.1)-(5.2)), it
suffices to show $\delta(x)^{1/2}\nabla^\ell (v\varphi-w)
\in L^2(\OO)$, where $\delta(x)=\text{dist}(x,\partial\OO)$.
But this is clear since $\delta(x)^{1/2}
\nabla^\ell(v\varphi)\in L^2(\OO)$
by the square function estimates as well as the assumption
$(\nabla^{\ell-1} v)^*\in L^2(\partial\OO)$, and $\nabla^\ell w\in L^2(\OO)$
by singular integral estimates \cite{St} and $\nabla^{\ell-1} v\in L^2(\OO)$.
Thus, by the $L^2$ uniqueness and estimate (1.6),
$$
\aligned
\int_{\Delta(Q_0,r)}
&|\Cal{M}_1(\nabla^{\ell-1}(v\varphi -w))|^p\, d\sigma
\le \int_{\partial\OO}
|(\nabla^{\ell-1}(v\varphi-w))^*|^p\, d\sigma\\
&\le C\int_{\partial\OO}
|\nabla^{\ell-1}(v\varphi-w)|^p\, d\sigma
= C\, \int_{\partial\OO} |\nabla^{\ell-1}w|^p\, d\sigma
\endaligned
\tag 4.18
$$
where we also used the fact $\nabla^{\ell-1}(v\varphi)=0$ on $\partial\OO$.

Let $p=q(n-1)/(n-q)$. Note that
$\frac{n-1}{p}=\frac{n}{q}-1$. It follows from (4.15) and Lemma 4.3 that
$$
\aligned
&\left(\frac{1}{r^{n-1}}\int_{\Delta(Q_0,5r)}
|\nabla^{\ell-1} w|^p\, d\sigma
\right)^{1/p}\\
&\ \ \ \ \ \ \le C\, r
\left(\frac{1}{r^n}
\int_{T(Q_0,3r)}
|\sum_{|\gamma|\le\ell-1}
r^{|\gamma|-\ell} \big|D^\gamma v|\big|^qdx\right)^{1/q}.
\endaligned
\tag 4.19
$$
Clearly in (4.19) we may replace $q$ by $\bar{q}=\max(q,2)$.
By Poincar\'e inequality (2.6), this gives
$$
\aligned
\left(\frac{1}{r^{n-1}}\int_{\Delta(Q_0,5r)}
|\nabla^{\ell-1} w|^p\, d\sigma\right)^{1/p}
&\le C\, \left(\frac{1}{r^n}
\int_{T(Q_0,3r)}
|\nabla^{\ell-1} v|^{\bar{q}}\, dx\right)^{1/\bar{q}}\\
&\le C\, \left(\frac{1}{r^{n-1}}
\int_{\Delta(Q_0,3r)}
|(\nabla^{\ell-1} v)^*|^{\bar{q}}\, d\sigma
\right)^{1/\bar{q}}.
\endaligned
\tag 4.20
$$

Finally, if $Q\in\partial\OO\setminus \Delta(Q_0,5r)$, we use
estimate (4.12) to obtain
$$
\aligned
|\nabla^{\ell-1} w(Q)|
& \le \frac{C}{|Q-Q_0|^{n+\ell-1}}
\int_{T(Q_0,3r)}
\sum_{|\gamma|\le \ell-1}r^{|\gamma|}
|D^\gamma v(y)|\, dy\\
&\le
\frac{C\, r^{n+\ell-1}}{|Q-Q_0|^{n+\ell-1}}
\left(\frac{1}{r^{n-1}}
\int_{\Delta(Q_0,3r)}
|(\nabla^{\ell-1} v)^*|^2\, d\sigma\right)^{1/2},
\endaligned
\tag 4.21
$$
as in (4.16).
It follows that
$$
\left(\frac{1}{r^{n-1}}
\int_{\partial\OO\setminus\Delta(Q_0,5r)}
|\nabla^{\ell-1} w|^p\, d\sigma\right)^{1/p}
\le C\, 
\left(\frac{1}{r^{n-1}}
\int_{\Delta(Q_0,3r)}
|(\nabla^{\ell-1} v)^*|^2\, d\sigma\right)^{1/2}.
\tag 4.22
$$
In view of (4.8), (4.14), (4.17), (4.20) and (4.22), we have proved that
$$
\left(\frac{1}{r^{n-1}}
\int_{\Delta(Q_0,r)}
|(\nabla^{\ell-1} v)^*|^p\, d\sigma\right)^{1/p}
\le C
\left(\frac{1}{r^{n-1}}
\int_{\Delta(Q_0,3r)}
|(\nabla^{\ell-1} v)^*|^{\bar{q}}\, d\sigma
\right)^{1/\bar{q}},
\tag 4.23
$$
where $\bar{q}=\max(q,2)$ and $\frac{n-1}{p}=\frac{n}{q}-1$.
Observe that for $p\ge 2$,
$$
\frac{1}{q}-\frac{1}{p}
=\frac{1}{n}(1-\frac{1}{p})\ge \frac{1}{2n}.
$$
With this, one may iterate estimate (4.23) to obtain
$$
\left(\frac{1}{r^{n-1}}
\int_{\Delta(Q_0,cr)}
|(\nabla^{\ell-1} v)^*|^p\, d\sigma\right)^{1/p}
\le C
\left(\frac{1}{r^{n-1}}
\int_{\Delta(Q_0,r)}
|(\nabla^{\ell-1} v)^*|^2\, d\sigma
\right)^{1/2},
\tag 4.24
$$
starting with $q=2$. This is possible since
 the $L^p$ solvability of the Dirichlet problem (1.4)
implies the $L^s$ solvability for any $2<s<p$.

By covering $\Delta(Q,r)$ with sufficiently
small surface balls $\{ \Delta(Q_j,cr)\}$,
it is easy to see that
estimate (4.24) is equivalent to the H\"older condition (1.8).
The proof is finished.
\enddemo

\remark{\bf Remark 4.25}
The proof of Theorem 4.1 would be much simpler if one assumes that
for the given $p$, the $L^p$ Dirichlet problem is uniquely solvable for
all Lipschitz domains. In this case, one may use the localization
techniques in \cite{DK1} and  apply the
$L^p$ estimate (1.6) on the domain
$\{ (x^\prime,x_n):\
|x^\prime|\le \rho r\ \text{ and }
\psi(x^\prime)<x_n<\psi(x^\prime)+ \rho r\}$
for $\rho\in (1,2)$.
\endremark

From Theorems 3.2 and 4.1 as well as the self-improving property
of the reverse H\"older condition (1.8), we may deduce the
following.

\proclaim{\bf Corollary 4.26}
Let $\CL(D)$ be an elliptic operator given by (1.1). Let $\OO$
be a bounded Lipschitz domain in $\br^n$, $n\ge 4$.
Then the set of exponents $p\in (2,\infty)$
for which the $L^p$ Dirichlet problem (1.4) on $\OO$
is uniquely solvable is an open interval
$(2,q)$ with $2<q\le \infty$.
\endproclaim

\medskip

\centerline{\bf 5. The Proof of Theorem 1.10}

In view of Theorem 3.2, to prove Theorem 1.10,
 it suffices to show that condition
(1.11) implies the reverse H\"older condition (1.8)
for $p$ in the range given by (1.12). 
To do this, we will use the regularity estimate (1.13).
The proof also depends on the following square
function estimates established in \cite{DKPV} for solutions
of $\CL(D)u=0$ in $\OO$,
$$
\align
\| S(\nabla^{\ell-1} u)\|_{L^p(\partial\OO)}
&\le C\, \| (\nabla^{\ell-1} u)^*\|_{L^p(\partial\OO)},\tag 5.1\\
\| (\nabla^{\ell-1} u)^*\|_{L^p(\partial\OO)}
&\le C\, 
\| S(\nabla^{\ell-1} u)\|_{L^p(\partial\OO)}
+C\, |\nabla ^{\ell -1} u (P_0)|\, |\partial\OO|^{1/p},\tag 5.2\\
\endalign
$$
where $0<p<\infty$, $P_0\in \OO$ and $C$ depends on
$n$, $m$, $\ell$, $\mu$, $P_0$ and the Lipschitz character of $\OO$. 

We first recall that
the square
function $S(w)$ is defined by
$$
S(w)(Q)=\left\{\int_{\Gamma(Q)}
\frac{|\nabla w (x)|^2}{|x-Q|^{n-2}}\, dx\right\}^{1/2}
\ \ \ \ \text{ for } Q\in \partial\OO.
\tag 5.3
$$
Let
$$
\widetilde{S} (w)(Q)
=\left\{ \int_{\Gamma(Q)}
\frac{|\nabla^2 w (x)|^2}{|x-Q|^{n-4}}\, dx
\right\}^{1/2}
\ \ \ \ \ \ \text{ for } 
Q\in\partial\OO.
\tag 5.4
$$
It follows from Lemma 2 on p.216 of \cite{St} that $\widetilde{S}
(w)(Q)$ is bounded
 by $C\, S(w)(Q)$ plus an interior term. Thus, by (5.2),
$$
\|
(\nabla^{\ell -1} u)^*\|_{L^p(\partial\OO)}
\le 
C\, 
\|\widetilde{S} (\nabla^{\ell -1} u)\|_{L^p(\partial\OO)}
+C\,\|
(\nabla ^{\ell -1} u)^*\|_{L^2(\partial\OO)}
|\partial\OO|^{\frac{1}{p}-\frac{1}{2}}.
\tag 5.5
$$

\proclaim{\bf Lemma 5.6} Let $p>2$. Then for any $\gamma\in (0,1)$
and $w\in C^2(\OO)$, we
have
$$
\int_{\partial\OO}
|\widetilde{S}(w)|^p\, d\sigma
\le C_\gamma
\left\{\text{diam} (\OO)\right\}^\gamma
\int_\OO |\nabla^2 w (x)|^p
\big[\delta(x)\big]^{2p-1-\gamma}\, dx,
\tag 5.7
$$
where $\delta(x)=\text{dist}(x,\partial\OO)$.
\endproclaim

\demo{Proof}
Write
$$
\widetilde{S}(w)(Q)
=\left\{
\int_{\Gamma(Q)}
\frac{|\nabla^2 w(x)|^2}{
|x-Q|^{\frac{2(n+\gamma)}{p}-4}}\cdot
\frac{dx}{|x-Q|^{\frac{(p-2)n-2\gamma}{p}}}\right\}^{1/2}.
\tag 5.8
$$
Using H\"older's inequality with exponents
$p/2$ and $(p/2)^\prime=\frac{p}{p-2}$ in (5.8), we obtain
$$
\widetilde{S}(w)(Q)
\le C\, \left\{ \text{diam}(\Omega)\right\}^{\frac{\gamma}{p}}
\left\{ \int_{\Gamma (Q)} \frac{|\nabla^2 w (x)|^p}
{|x-Q|^{n+\gamma -2p}}\, dx\right\}^{1/p}.
\tag 5.9
$$
From this, inequality (5.7) follows easily by integration.
\enddemo

\proclaim{\bf Lemma 5.10}
Let $p>2$. Suppose $\CL(D)u=0$ in $\OO$.
Then for any $\gamma\in (0,1)$,
$$
\aligned
&\int_{\partial\OO} |(\nabla^{\ell-1} u)^*|^p \, d\sigma
\le C\, \left\{ \int_{\partial\OO}
|(\nabla^{\ell -1} u)^*|^2\, d\sigma\right\}^{p/2}\, |\partial\OO|^{1-\frac{
p}{2}}\\
&+C_\gamma
\left\{ \text{diam}(\OO)\right\}^\gamma
\sup_{x\in\OO}
|\nabla^{\ell+1}u(x)|^{p-2}
\big[\delta(x)\big]^{2p-2-\gamma}
\int_{\partial\OO}
|(\nabla^\ell u)^*|^2\, d\sigma.
\endaligned
$$
\endproclaim 

\demo{Proof}
It follows from (5.7) that
$$
\aligned
\int_{\partial\OO}
&|\widetilde{S}(\nabla^{\ell -1} u)|^p\, d\sigma
\le
C_\gamma \left\{ \text{diam}(\OO)\right\}^\gamma
\int_\OO 
|\nabla^{\ell +1} u|^p
\big[\delta(x)\big]^{2p-1-\gamma}\, dx\\
&\le C_\gamma
\left\{ \text{diam}(\OO)\right\}^\gamma
\sup_{x\in\OO}
|\nabla^{\ell +1} u(x)|^{p-2}
\big[\delta(x)\big]^{2p-2-\gamma}
\int_\OO |\nabla^{\ell +1} u|^2
\,\delta(x)\, dx\\
&\le
C_\gamma\left\{\text{diam}(\OO)\right\}^\gamma
\sup_{x\in\OO}
|\nabla^{\ell+1}u(x)|^{p-2} \big[\delta(x)\big]^{2p-2 -\gamma}
\int_{\partial\OO}
|(\nabla^\ell u)^*|^2\, d\sigma,
\endaligned
$$
where we have used square function estimate
 (5.1) (with $p=2$) in the last inequality.
This, together with (5.5), gives the desired
estimate in Lemma 5.10.
\enddemo
 
We are now ready to give the proof of Theorem 1.10.
\demo{\bf Proof of Theorem 1.10}
We begin by fixing $\Delta(Q_0,r)$ with $Q_0\in \partial\OO$ and $0<r<r_0$.
 Let $v$ be a solution
of $\CL(D)v=0$ in $\OO$ with the properties
$(\nabla^{\ell-1} v)^*\in L^2(\partial\OO)$ and 
$D^\alpha v=0$ on $\Delta(Q_0,3r)$ for $|\alpha|\le \ell-1$.
Note that by condition (1.11) and interior estimate (2.2),
 for any $x\in T(Q_0,r)$,
$$
\big[\delta(x)\big]^2|\nabla^{\ell+1} v(x)|\le C\, 
\left(\frac{\delta(x)}{r}\right)^{\frac{\lambda-n}{2}}
\left(\frac{1}{r^n}\int_{T(Q_0,2r)}
|\nabla^{\ell-1} v|^2\, dx\right)^{1/2}.
$$
It then follows that for any $x\in T(Q_0,r)$,
$$
|\nabla^{\ell+1} v(x)|
\le \frac{C}{\big[\delta(x)\big]^2}
\left(\frac{\delta(x)}{r}\right)^{\frac{\lambda-n}{2}}
\left(\frac{1}{r^{n-1}}
\int_{\Delta(Q_0,2r)}
|(\nabla^{\ell-1}v)^*|^2\, d\sigma\right)^{1/2}.
\tag 5.11
$$

By rotation and translation,
we may assume that $Q_0=0$ and $r_0=r_0(n,\OO)>0$ is so small that
$$
\aligned
B(0,C_0\,r_0)\cap \OO
&=B(0,C_0\,r_0)\cap \big\{(x^\prime,x_n)\in \br^{n-1}\times\br:\ 
x_n>\psi(x^\prime)\big\},\\
B(0,C_0r_0)\cap \partial\OO
&=B(0,C_0\,r_0)\cap \big\{(x^\prime,\psi(x^\prime)):\ 
x^\prime\in\br^{n-1}\big\}
\endaligned
$$
where $\psi$ is a Lipschitz function on $\br^{n-1}$.
For $\rho\in (1,4)$, with slightly abused notation, we let
$$
\aligned
I_{\rho r}&=\big\{ (x^\prime,\psi(x^\prime)):\ |x^\prime|<\rho c_2 r\big\},\\
Z_{\rho r} &=\big\{ (x^\prime, x_n):\
|x^\prime|<\rho c_2 r\ \text{ and }\
\psi(x^\prime)<x_n<\psi(x^\prime) +\rho c_2 r\big\},
\endaligned
\tag 5.12
$$
where $c_2=c_2(n,\| \nabla\psi\|_\infty)>0$ is small so that
$I_{3r}\subset \Delta(0,r)$ and $Z_{3r}\subset B(0,r)\cap \OO$.
Let $\Cal{M}_1$ and $\Cal{M}_2$ be the operators defined
by (4.6). As in (4.7) and (4.8), it is easy to see that
$$
\left(\frac{1}{r^{n-1}}
\int_{I_r} |\Cal{M}_2(\nabla^{\ell-1} v)|^p\, d\sigma\right)^{1/p}
\le C\, \left(\frac{1}{r^{n-1}}
\int_{\Delta(0,2r)} |(\nabla^{\ell-1} v)^*|^2\,
d\sigma\right)^{1/2},
\tag 5.13
$$
by interior estimates.

To estimate $\Cal{M}_1(\nabla^{\ell-1} v)$ on $I_r$, we apply
Lemma 5.10 to solution 
$v$ on the Lipschitz domain $Z_{\rho r}$ for $\rho\in (3/2,2)$.
This gives
$$
\aligned
&\frac{1}{r^{n-1}}\int_{I_r} |\Cal{M}_1(\nabla^{\ell-1} v)|^p\, d\sigma
\le \frac{1}{r^{n-1}}
\int_{\partial Z_{\rho r}}
|(\nabla^{\ell-1} v)^*_\rho|^p\, d\sigma\\
&\le C\,\left(\frac{1}{r^{n-1}}
\int_{\partial Z_{\rho r}}
|(\nabla^{\ell-1} v)^*_\rho|^2\, d\sigma\right)^{p/2}\\
&+C_\gamma\, r^{\gamma}
\sup_{x\in Z_{\rho r}} |\nabla^{\ell+1}v(x)|^{p-2}\big[\delta_\rho (x)\big]^{
2p-2-\gamma}
\cdot\frac{1}{r^{n-1}}
\int_{\partial Z_{\rho r}}
|(\nabla^\ell v)^*_\rho|^2\, d\sigma,
\endaligned
\tag 5.14
$$
where $\delta_\rho (x)=\text{dist}(x,\partial Z_{\rho r})$ and
$(\nabla^\ell v)_\rho^*$ denotes the non-tangential maximal
function of $\nabla^\ell v$ with respect to the domain $Z_{\rho r}$.
By regularity estimate (1.13),
$$
\int_{\partial Z_{\rho r}} |(\nabla^{\ell} v)^*_\rho|^2\, d\sigma
\le C\, \int_{\partial Z_{\rho r}}
|\nabla_t\nabla^{\ell-1} v|^2\, d\sigma
\le C\, \int_{\OO\cap \partial Z_{\rho r}}|\nabla^\ell v|^2\, d\sigma,
\tag 5.15
$$
since $\nabla^{\ell-1} v =0$ on
$\Delta(0,3r)$.

Note that $\delta_\rho (x)\le \delta(x)\le C\, r$ for $x\in Z_{\rho r}$.
Also observe that the condition (1.12) for $p$ is equivalent to
$$
\frac{\lambda -n}{2}\cdot (p-2) +2>0.
$$
Thus we may choose $\gamma>0$ so small that
$$
\frac{\lambda -n}{2}\cdot (p-2) +2-\gamma>0.
$$
By (5.11), this implies that
$$
\aligned
r^\gamma \sup_{x\in Z_{\rho r}}
|\nabla^{\ell+1} v(x)|^{p-2}
&\big[\delta_\rho(x)\big]^{2p-2-\gamma}\\
&\le C\, r^2
\left(\frac{1}{r^{n-1}}
\int_{\Delta(0,2r)}
|(\nabla^{\ell-1} v)^*|^2\, d\sigma\right)^{\frac{p-2}{2}}.
\endaligned
\tag 5.16
$$
In view of (5.14), (5.15) and (5.16), we have proved that
$$
\aligned
&\left(\frac{1}{r^{n-1}}
\int_{I_r}
|\Cal{M}_1(\nabla^{\ell-1} v)|^p\, d\sigma \right)^{2/p}
\le \frac{C}{r^{n-1}}
\int_{\OO\cap \partial Z_{\rho r}}|\nabla^{\ell-1} v|^2\, d\sigma\\
&\ \ \ \
+C\, r^{\frac{4}{p}}
\left(\frac{1}{r^{n-1}}
\int_{\Delta(0,2r)}
|(\nabla^{\ell-1} v)^*|^2\, d\sigma\right)^{\frac{p-2}{p}}
\left(\frac{1}{r^{n-1}}
\int_{\OO\cap \partial Z_{\rho r}}
|\nabla^\ell v|^2\, d\sigma\right)^{2/p}.
\endaligned
$$
By H\"older's inequality, it follows that
$$
\aligned
&\left(\frac{1}{r^{n-1}}\int_{I_r}
|\Cal{M}_1(\nabla^{\ell-1} v)|^p\, d\sigma\right)^{2/p}
\le \frac{C}{r^{n-1}}\int_{\Delta(0,2r)}
|(\nabla^{\ell-1} v)^*|^2\, d\sigma\\
&\ \ \ \ \ \ \ \ \ \ \ \ \ \ \ \ \ \ \ \ \
+\frac{C}{r^{n-1}}\int_{\OO\cap \partial Z_{\rho r}}
|\nabla^{\ell-1} v|^2\, d\sigma
+\frac{C}{r^{n-3}}
\int_{\OO\cap\partial Z_{\rho r}}
|\nabla^\ell v|^2\, d\sigma.
\endaligned
$$
Integrating the above inequality in $\rho\in (3/2,2)$, we obtain
$$
\aligned
&\left(\frac{1}{r^{n-1}}\int_{I_r}
|\Cal{M}_1(\nabla^{\ell-1} v)|^p\, d\sigma\right)^{2/p}
\le \frac{C}{r^{n-1}}\int_{\Delta(0,2r)}
|(\nabla^{\ell-1} v)^*|^2\, d\sigma\\
&\ \ \ \ \ \ \ \ \ \ \ \ \ \ \ \ \ \ \ \ \ \ \ \
+\frac{C}{r^n}\int_{Z_{2 r}}
|\nabla^{\ell-1} v|^2\, d\sigma
+\frac{C}{r^{n-2}}
\int_{Z_{2r}}
|\nabla^\ell v|^2\, d\sigma.
\endaligned
\tag 5.17
$$
Using Cacciopoli's inequality (2.10), it is easy to see that
the last two terms in the right side of (5.17) is dominated
by the first term. This, together with (5.13), gives
$$
\left(\frac{1}{I_r}
\int_{I_r} |(\nabla^{\ell-1} v)^*|^p\, d\sigma \right)^{1/p}
\le C\,
\left(\frac{1}{r^{n-1}}
\int_{\Delta(0,2r)}
|(\nabla^{\ell-1} v)^*|^2\, d\sigma\right)^{1/2}.
\tag 5.18
$$
By a simple covering argument, inequality (5.18) implies
the reverse H\"older condition (1.8).
The proof is complete.
\enddemo

\demo{\bf Proof of Corollary 1.14}
We will show that for any elliptic operator $\CL(D)$ given by (1.1),
condition (1.11) holds for some $\lambda>3$.
To this end,
let $v$ be a solution of $\CL(D) v=0$ in $\OO$ with the properties
$(\nabla^{\ell-1} v)^*\in L^2(\partial\OO)$ and 
$D^\alpha v=0$ on $\Delta(Q_0,5R)$ for $|\alpha|\le \ell-1$.
We will assume $Q_0=0$ and use the same notation as in the proof
of Theorem 1.10.

Let $0<r<R/2$. Note that by H\"older's inequality,
$$
\aligned
\int_{Z_r} |\nabla^{\ell-1}v|^2\, dx
&\le C\, r^3 \int_{I_r}
|\Cal{M}_1 (\nabla^\ell v)|^2 \, d\sigma\\
&\le C\, r^{3+(n-1)(1-\frac{2}{q})}
\left(\int_{I_{\rho R}}
|\Cal{M}_1(\nabla^\ell v)|^{q}\, d\sigma\right)^{2/q},
\endaligned
$$
where $\rho\in (1/2,1)$ and $q>2$. Choose $q>2$ so that
the regularity estimate (1.13) holds on the Lipschitz domain
$Z_{\rho R}$ uniformly for $\rho\in (1/2,1)$. 
It follows that
$$
\left(\int_{Z_r} |\nabla^{\ell-1}v|^2\, dx\right)^{q/2}
\le C\, r^{\frac{3q}{2}+(n-1)(\frac{q}{2}-1)}
\int_{\OO\cap \partial Z_{\rho R}}
|\nabla^\ell v|^{q}\, d\sigma.
\tag 5.19
$$
Integrating both sides of (5.19) in $\rho\in (1/2,1)$, we obtain
$$
\left(\int_{Z_r} |\nabla^{\ell-1}v|^2\, dx\right)^{q/2}
\le C\, r^{\frac{3q}{2}+(n-1)(\frac{q}{2}-1)}
\cdot\frac{1}{R}\int_{Z_R}
|\nabla^\ell v|^{q}\, dx.
$$
We may assume that inequality (2.12) holds for this $q$. Hence,
$$
\aligned
\int_{Z_r} |\nabla^{\ell-1}v|^2\, dx
&\le C\, r^{n+2}
\cdot
\left(\frac{r}{R}\right)^{(1-n)(\frac{2}{q})}
\cdot\frac{1}{R^n}
\int_{Z_{2R}}
|\nabla^\ell v|^2\, dx\\
& \le C \left(\frac{r}{R}\right)^{3 +(n-1)(1-\frac{2}{q})}
\int_{Z_{3R}} |\nabla^{\ell-1} v|^2\, dx,
\endaligned
\tag 5.20
$$
where we have used the Cacciopoli's inequality (2.10) in the last inequality.
By a covering argument, estimate (5.20) implies the condition (1.11)
with $\lambda=3+ (n-1)(1-\frac{2}{q})>3$.
\enddemo

\remark{\bf Remark 5.21}
It would be very interesting to know whether condition (1.11) in Theorem 1.10
is also necessary in the case $p>2(n-1)/(n-3)$.
That is, 
does  the $L^p$ solvability of the Dirichlet problem on $\OO$
implies  condition (1.11) for all 
$\lambda< n-\frac{4}{p-2}$? We point out that
it is not hard to see that the $L^p$ solvability implies condition (1.11)
for $\lambda=n-\frac{2(n-1)}{p}$.
However, this is a weaker statement if $p>2(n-1)/(n-3)$.
\endremark

\enddocument

\end